\documentclass[12pt,a4paper]{article}

\usepackage{a4wide,amssymb,amsmath,amsthm,xspace,epsfig, amsfonts}

\usepackage[left=1.5cm, right=1.5cm, top=1.5cm]{geometry} 

\usepackage{tikz-cd}

\begin{document}

\newcommand{\N}{\mathbb{N}}
\newcommand{\R}{\mathbb{R}}
\newcommand{\Z}{\mathbb{Z}}
\newcommand{\Q}{\mathbb{Q}}
\newcommand{\C}{\mathbb{C}}
\newcommand{\PP}{\mathbb{P}}

\newcommand{\LL}{\Bbb L}
\newcommand{\OO}{\mathcal{O}}
\newcommand{\DD}{\mathcal{D}}

\newcommand{\esp}{\vskip .3cm \noindent}
\mathchardef\flat="115B

\newcommand{\lev}{\text{\rm Lev}}

\def\ut#1{$\underline{\text{#1}}$}
\def\CC#1{${\cal C}^{#1}$}
\def\h#1{\hat #1}
\def\t#1{\tilde #1}
\def\wt#1{\widetilde{#1}}
\def\wh#1{\widehat{#1}}
\def\wb#1{\overline{#1}}

\def\restrict#1{\bigr|_{#1}}

\def\ufin#1#2{\mathsf{U_{fin}}\bigl({#1},{#2}\bigr)}
\def\sfin#1#2{\mathsf{S_{fin}}\bigl({#1},{#2}\bigr)}
\def\sone#1#2{\mathsf{S_{1}}\bigl({#1},{#2}\bigr)}

\def\ufinw#1#2{\mathsf{U_{fin}^w}\bigl({#1},{#2}\bigr)}
\def\sfinw#1#2{\mathsf{S_{fin}^w}\bigl({#1},{#2}\bigr)}
\def\sonew#1#2{\mathsf{S_{1}^w}\bigl({#1},{#2}\bigr)}

\def\ch#1#2{\left(\begin{array}{c}#1 \\ #2 \end{array}\right)}

\newtheorem{lemma}{Lemma}[section]

\newtheorem{thm}[lemma]{Theorem}

\newtheorem{defi}[lemma]{Definition}
\newtheorem{conj}[lemma]{Conjecture}
\newtheorem{cor}[lemma]{Corollary}
\newtheorem{prop}[lemma]{Proposition}
\newtheorem*{prob}{Problem}
\newtheorem{q}[lemma]{Question}
\newtheorem*{rem}{Remark}
\newtheorem{examples}[lemma]{Examples}
\newtheorem{example}[lemma]{Example}

\title{On (non-Menger) spaces whose closed nowhere dense subsets are Menger}

\date{\empty}
\author{Mathieu Baillif and Santi Spadaro}
\maketitle

\abstract{\footnotesize
          A space $X$ is od-Menger if it satisfies $\ufin{\Delta_X}{\OO_X}$, where $\OO_X,\Delta_X$
          are the collection of covers of $X$ by respectively open subsets and open dense subsets.
          We show that under {\bf CH}, there is a refinement of the usual topology on a subset of the reals which yields a 
          hereditarily Lindel\"of, od-Menger, non-Menger, $0$-dimensional, first countable space.
          We also investigate the properties of spaces which are od-Menger but not Menger.}

\section{Introduction}

We follow the set theoretic convention of denoting the integers by $\omega$.
In this paper, by {\em space} we mean a $T_1$ topological space, and a {\em cover} of a space 
is a collection of open sets whose union contains the space. 
$\OO_X,\Delta_X$ respectively mean the collection of covers and od-covers of the topological space $X$, where
a cover is an od-cover iff every member is open dense.
We say that a subset $\mathcal{U}$ 
of the subsets of $X$
is a weak cover of $X$ if $\cup\mathcal{U}$ is dense in $X$, and let
$\DD_X$ be the collection of weak covers of $X$.
Let $\mathcal{A},\mathcal{B}$ be collections of subsets 
of $X$. We recall the definitions of the following classical properties:
\begin{itemize}
  \item[\ ]
  $\ufin{\mathcal{A}}{\mathcal{B}}$: 
  For each sequence $\langle \mathcal{U}_n \colon n\in\omega\rangle$ of members of $\mathcal{A}$, 
  there are finite $\mathcal{F}_n\subset \mathcal{U}_n$ such that $\{\cup \mathcal{F}_n\,:\,n\in\omega\} \in \mathcal{B}$.
  \item[\ ]
  $\sone{\mathcal{A}}{\mathcal{B}}$: 
  For each sequence $\langle \mathcal{U}_n \colon n\in\omega\rangle$ of members of $\mathcal{A}$ 
  there is $U_n\in\mathcal{U}_n$ such that $\{U_n\,:\,n\in\omega\} \in \mathcal{B}$.
\end{itemize}

A space $X$ satisfying $\ufin{\OO_X}{\OO_X}$~(resp. $\ufin{\OO_X}{\DD_X}$) [resp. $\ufin{\Delta_X}{\OO_X}$]
is called {\em Menger} (resp. {\em weakly Menger}) [resp. od-Menger],
and one satisfying $\sone{\OO_X}{\OO_X}$~(resp. $\sone{\OO_X}{\DD_X}$) [resp. $\sone{\Delta_X}{\OO_X}$]
is called {\em Rothberger} (resp. {\em weakly Rothberger}) [resp. od-Rothberger].
If $X$ is clear from the context, we may omit the subscript on $\OO_X,\Delta_X,\DD_X$
and simply say ``Let $X$ be such that $\ufin{\OO}{\OO}$ holds'', for instance.
Notice in passing that $\sone{\Delta_X}{\DD_X}$ holds for any space $X$.
The following lemma follows immediately from the definitions.
Recall that a space is {\em od-Lindel\"of} iff any cover by open dense sets has a countable subcover (or equivalently
if any nowhere dense subset is Lindel\"of), and {\em weakly Lindel\"of} iff any cover countains a countable subfamily which is 
a weak cover.
\begin{lemma} 
   Let ``has dM'' and ``has dR'' stand for ``has a dense Menger subspace'' and 
   ``has a dense Rothberger subspace'', respectively.
   The implications below hold for 
   any space $X$. 
   \begin{center}
   \begin{tikzcd}
   & & 
   X \text{ is separable} 
       \arrow[d, Rightarrow]
   &
   \\
   \sone{\Delta_X}{\OO_X} 
     \arrow[r, Leftarrow]
     \arrow[d, Rightarrow]
   &
   \sone{\OO_X}{\OO_X} 
    \arrow[r, Rightarrow] 
    \arrow[d, Rightarrow]
   & 
   X\text{ has dR} 
     \arrow[r, Rightarrow]
     \arrow[d, Rightarrow]
   &
   \sone{\OO_X}{\DD_X}
     \arrow[d, Rightarrow]
   \\
   \ufin{\Delta_X}{\OO_X} 
     \arrow[r, Leftarrow] 
     \arrow[d, Rightarrow]
   &
   \ufin{\OO_X}{\OO_X} 
     \arrow[r, Rightarrow] 
     \arrow[d, Rightarrow]
   & 
   X\text{ has dM} 
     \arrow[r, Rightarrow]
   &
   \ufin{\OO_X}{\DD_X}
     \arrow[d, Rightarrow]
   \\
   X\text{ is od-Lindel\"of} 
      \arrow[r, Leftarrow]
   &
   X\text{ is Lindel\"of} 
     \arrow[rr, Rightarrow]
   &
   &
   X\text{ is weakly Lindel\"of}  
   \end{tikzcd}
   \end{center}
\end{lemma} 
There are obvious examples of od-Rothberger non-Menger spaces: uncountable discrete spaces, for instance.
However, the question is more subtle for Lindel\"of spaces, and
the first author asked in \cite{meszigues-od-sel} whether a Lindel\"of od-Menger space is actually Menger.
We were unaware that a space of Sakai
\cite[Thm 2.4 (b)]{Sakai:2014}, built with the continuum hypothesis {\bf CH} 
for other (though similar) purposes, could very easily be shown to be a counter-example
(see below after Corollary \ref{ex:main} for more details). 
In this note, we present another {\bf CH} counter-example, due to the second author,
and prove some properties of spaces that are od-Menger and non-Menger.
First, we notice that the
question of whether an od-Lindel\"of space is actually Lindel\"of was settled in 
the next easy\footnote{Easy in hindsight,
because Theorem \ref{thm:MWB} had escaped the first author who spent an entire paper trying to prove it, only
to find convoluted proofs of weaker results.} theorem, which is proved in 
\cite{meszigues-od-sel, Blair:1983, MillsWattel}.

\begin{thm}[{Mills, Wattel and Blair}]\label{thm:MWB}
   If $X$ is an od-Lindel\"of $T_1$ space, then $X$ either has an uncountable clopen discrete subset or $X$ is Lindel\"of.
   In particular, the subspace of non-isolated points is Lindel\"of.
\end{thm}

We now gather some lemmas (classical or easy) that will be of use later.
The first one is almost immediate, 
the easy proof for $\mathsf{U_{fin}}$ given in \cite[Lemma 3.6]{meszigues-od-sel} works as well for 
$\mathsf{S_{1}}$.

\begin{lemma} \label{lemmaequiv}
  Let $X$ be a space and $\mathsf{P}\in\{\mathsf{U_{fin}},
  \mathsf{S_{1}}\}$. The items below are equivalent:\\
  (a) $X$ satisfies $\mathsf{P}(\Delta,\OO)$,\\
  (b) $X$ satisfies $\mathsf{P}(\Delta,\Delta)$,\\
  (c) $X$ satisfies $\mathsf{P}(\Delta_1,\OO)$, where $\Delta_1$ is the collection of open covers with at least one dense member,\\
  (d) any closed subset of $X$ satisfies $\mathsf{P}(\Delta,\OO)$,\\
  (e) any closed nowhere dense subset of $X$ satisfies $\mathsf{P}(\OO,\OO)$.\\
\end{lemma}

The next lemma follows immediately.
\begin{lemma}\label{lemma:DeltaOandOD}
   For any space $X$ and $\mathsf{P}\in\{\mathsf{U_{fin}},
   \mathsf{S_{1}}\}$,
   $$ \mathsf{P}(\Delta_X,\OO_X) \text{ and } \mathsf{P}(\OO_X,\DD_X) \,\Longrightarrow\, \mathsf{P}(\OO_X,\OO_X) $$
\end{lemma}
\proof
   Given a sequence of open covers, extract a weak cover from the odd indices and
   from the even indices a cover of the closed nowhere dense remaining set.
\endproof

Recall the partial order of eventual dominance in $^\omega\omega$: $f\le^* g$ iff there is some $m\in\omega$
such that $f(n)\le g(n)$ whenever $n\ge m$.
A family $A\subset\,^\omega\omega$ is {\em unbounded} if
there is no $f\in\,^\omega\omega$ with $g\le^* f$ for each $g\in A$, and {\em dominating} iff
for each $g\in\,^\omega\omega$, there is $f\in A$ with $g\le^* f$.
Then $\mathfrak{b}$ and $\mathfrak{d}$ are respectively 
the smallest cardinal of an unbounded and of a dominating family in $^\omega\omega$,
while $\mathfrak{c}$ denotes the cardinality of the continuum.
It is classical that $\omega_1\le\mathfrak{b}\le\mathfrak{d}\le\mathfrak{c}$. It is also
well known that $\mathfrak{b}$ is a regular cardinal.

\begin{thm}[Folklore, probably]
   \label{thm:led}
   The following hold. \\
   (a) A Lindel\"of space of cardinality $<\mathfrak{d}$ is Menger.\\
   (b) A Lindel\"of space which is
   the union of $<\mathfrak{b}$ Menger subspaces is Menger.
\end{thm}

Finally, a last classical observation given without proof.
\begin{lemma}\label{lemma:classicalequiv}
   Let $\mathcal{B}\subset\OO_X$ and $\mathcal{A}$ be a collection of subsets of $X$. Then
   $X$ satisfies 
   $\ufin{\mathcal{A}}{\mathcal{B}}$ (resp. $\sone{\mathcal{A}}{\mathcal{B}}$) iff
   for each sequence $\langle \mathcal{U}_n \colon n\in\omega\rangle$ of members of $\mathcal{A}$
   there are finite (resp. one element subsets) $\mathcal{F}_n\subset \mathcal{U}_n$ such that for
   each $x\in X$ there are infinitely many $n$ such that $x\in \cup \mathcal{F}_n$. 
\end{lemma}


\section{Properties of spaces that are od-Menger but non-Menger}

Recall that a space $X$ has a $G_\delta$ diagonal iff there is a sequence $\langle \mathcal{D}_n\,:\,n\in\omega\rangle$
of open covers of $X$ such that for all $x\in X$, $\displaystyle\cap_{n\in\omega} St(x,\mathcal{D}_n) = \{x\}$
(where $St(x,\mathcal{D}_n)$ is the union of members of $\mathcal{D}_n$ containing $x$),
or equivalently iff the diagonal is a $G_\delta$ subset of $X^2$.
The sequence $\langle \mathcal{D}_n\,:\,n\in\omega\rangle$ is then called a $G_\delta$-diagonal sequence.

\begin{prop}\label{prop:real_example}
   If there is a Lindel\"of $T_1$ space with a $G_\delta$ diagonal which satisfies 
   $\ufin{\Delta}{\OO}$ but not $\ufin{\OO}{\OO}$,
   then there is one which is a dominating subset of $^\omega\omega$
   endowed with a topology finer than the Euclidean one.
\end{prop}

\proof
The proof is similar to many arguments about (non-)Menger spaces.
We will consider finite and infinite sequences in $^{<\omega}\omega$ and $^\omega \omega$,
for clarity we will put a bar above infinite ones, 
writing for instance $\wb{\sigma}\in\,^\omega \omega$ and $\tau\in\,^{<\omega}\omega$.
Let $\mathfrak{U} = \langle \mathcal{U}_n \colon n\in\omega\rangle$ be a sequence of open covers of $X$
witnessing that $X$ does not satisfy $\ufin{\OO}{\OO}$.
We may assume that $\mathcal{U}_n=\{U_{n,k}\,:\,k\in\omega\}$ for each $n$.
It will be convenient to 
consider refinements of (a subcover of) each $\mathcal{U}_n$ defined as follows.
Let $\langle \mathcal{D}_n\,:\,n\in\omega\rangle$ be a $G_\delta$-diagonal sequence.
We may assume that $\mathcal{D}_n$ is
countable for each $n$. 
Up to taking all intersections of members of $\mathcal{U}_n$ and $\mathcal{D}_n$ (and re-enumerate $\mathcal{U}_n$),
we may also
assume that, for each $n,k$, there is some $D\in\mathcal{D}_n$ such that $U_{n,k}\subset D$.
Given $\sigma\in\,^{<\omega}\omega$, set $U_\sigma = \cap_{k\in\text{dom}(\sigma)}U_{k,\sigma(k)}$ (which may be empty),
then $\mathcal{V}_n = \{ U_\sigma\,:\, \text{dom}(\sigma)=n+1\}$ is a cover of $X$ and a refinement of $\mathcal{U}_n$. 
Set $\mathfrak{V}=\langle \mathcal{V}_n \colon n\in\omega\rangle$.
\\
We fix a sequence $\mathfrak{W}=\langle\mathcal{W}_n\,:\,n\in\omega\rangle$  of families $\mathcal{W}_n= \{ W_\sigma\colon\sigma\in\,^{n+1}\omega\}$ of
pairwise disjoints open sets such that $\cup\mathcal{W}_n$ is dense in $X$ for each $n$ 
and $W_\sigma\subset U_\sigma \cap W_{\sigma\upharpoonright k+1}$
for each $k<n$. We allow $W_\sigma$ to be empty. 
This may be done by induction on $n$ (and using Zorn's Lemma) by considering maximal families 
of $2$-by-$2$ disjoint open sets such that each is contained in some $U_\sigma \cap W_{\sigma\upharpoonright n}$,
starting with $W_\varnothing = X$.
Then $\cup_n (X-\cup\mathcal{W}_n) = X - \cap_n \cup\mathcal{W}_n$ is a countable union of closed nowhere dense subspaces and
thus satisfies $\ufin{\OO}{\OO}$. We cover it with the union of finite subfamilies of each $\mathcal{V}_n$.
Let $Y$ be the closed space not covered.
Notice that each $\mathcal{W}_n$ is a cover of $Y$.
If there are finite subfamilies of each $\mathcal{W}_n$ whose union covers $Y$, then the same holds for $\mathcal{V}_n$
(and a fortiori for $\mathcal{U}_n$).
We can thus assume that 
\begin{equation}\label{eq:finitesub}\tag{*}
   \text{If }\mathcal{F}_n\subset\mathcal{W}_n\text{ are finite subfamilies, then }\cup_{n\in\omega}\cup\mathcal{F}_n \not\supset Y.
\end{equation}
Given $\wb{\sigma}\in\,^\omega\omega$, set $W_{\wb{\sigma}} = Y\bigcap\cap_n W_{\wb{\sigma}\upharpoonright n}$.
If non-empty, $W_{\wb{\sigma}}$ contains exactly one point $x_{\wb{\sigma}}$, 
because $W_{\wb{\sigma}\upharpoonright n+1}$ is a subset
of some $D\in\mathcal{D}_n$ and $\langle \mathcal{D}_n\,:\,n\in\omega\rangle$ is a $G_\delta$-diagonal sequence.
Moreover, by disjointness of $\mathcal{W}_n$,
each point $y\in Y$ is equal to a unique $x_{\wb{\sigma}}$.\\
By construction, the map 
$\mu:Y\to\, ^\omega\omega$ defined by 
$x_{\wb{\sigma}}\mapsto\wb{\sigma}$ is continuous and $1$-to-$1$. 
Moreover, the image of the open sets $W_\sigma$ are the basic open sets $[\sigma]$ of $^\omega\omega$.
It follows from (\ref{eq:finitesub})
that $\mu(Y)$ does not satisfy $\ufin{\OO}{\OO}$ and hence is a dominating subset of $^\omega\omega$.
\endproof

\begin{defi}\ \\
    $\bullet$ 
    If $X$ is a space, let $M(X)$ and $R(X)$ be its open subspaces 
    \begin{align*}
        M(X) &= \{x\in X\,:\, \exists U\ni x,\,U\text{ open and }\wb{U}\text{ Menger }\};\\
        R(X) &= \{x\in X\,:\, \exists U\ni x,\,U\text{ open and }\wb{U}\text{ Rothberger }\}.
    \end{align*}
    $\bullet$
    A space $X$ is everywhere non-Menger iff $\wb{U}$ is non-Menger for every 
    non-empty open $U\subset X$, that is,
    if $M(X) = \varnothing$; and $X$ is everywhere non-Rothberger iff $\wb{U}$ is non-Rothberger
    for every non-empty open $U\subset X$, that is,
    if $R(X) = \varnothing$.
\end{defi}

\begin{lemma}\label{lemma:odMdenseM}
   Let $X$ be a Lindel\"of space. \\
   (1)
   If $X$ is od-Menger and $M(X)$ is dense in $X$, then $X$ is Menger.\\
   (2)
   If $X$ is od-Rothberger and $R(X)$ is dense in $X$, then $X$ is Rothberger.
\end{lemma}
\begin{proof}
   \ \\
   (1)
   Let $X$ satisfy $\ufin{\Delta}{\OO}$.
   Let $\mathcal{W}$ be a cover of $M(X)$
   given by choosing $U_x\ni x$ with a Menger closure for each $x\in M(X)$.
   Let $\mathcal{V}$ be a maximal family of disjoint open sets refining $\mathcal{W}$.
   Then by maximality $\cup\mathcal{V}$ is dense in $M(X)$ and hence in $X$.
   Given any sequence $\mathcal{U}_n$ of covers of $X$, since $X$ is od-Menger we may cover 
   $X-\cup\mathcal{V}$ with the union of finite families of each $\mathcal{U}_n$ for even $n$.
   Let $D$ be the remaining uncovered closed subspace of $X$.
   Then $D\subset\cup\mathcal{V}$. By Lindel\"ofness of $D$ and disjointness of $\mathcal{V}$,
   $D$ is contained in a countable union of Menger subspaces of $X$, and hence is Menger and may be covered
   by the union of finite families of each $\mathcal{U}_n$ for odd $n$. This shows that $X$ is Menger. \\
   (2) Same proof.
\end{proof}

\begin{cor}
   \
   \\
   (1)
   If there is a Lindel\"of od-Menger space which is not Menger, then there is one
   which is everywhere non-Menger.\\
   (2) If there is a Lindel\"of od-Rothberger space which is not Menger, then there is one
   which is everywhere non-Menger.\\
   (3) If there is a Lindel\"of od-Rothberger space which is not Rothberger, then there is one
   which is everywhere non-Rothberger.
\end{cor}

\begin{proof} \
   \\
   (1)
   By Lemma \ref{lemma:odMdenseM},
   $M(X)$ is not dense in $X$, let $Y = \wb{\text{int} \left(X-M(X)\right) }$.
   Then $Y$ is Lindel\"of, od-Menger and non Menger, and if $U$ is open and $U\cap Y\not=\varnothing$, then 
   $\wb{U}\cap Y$ is not Menger.\\
   (2) Follows immediately by (1) since od-Rothberger implies od-Menger.\\
   (3) Same proof as (1), using $R(X)$ instead of $M(X)$.
\end{proof}

The next lemma is inspired by \cite[Lemma 13]{Tall:2011}.
\begin{lemma}
   \label{lemma:compact}
   Let $X$ be Lindel\"of and od-Menger, and $Y\subset X$ be the union of $<\mathfrak{d}$ compact subsets.
   Then $\wb{Y}$ is Menger. 
\end{lemma}
\begin{proof}
   $\wb{Y}$ is Lindel\"of, let $\mathcal{U}_n = \{U_{n,m}\,:\,m\in\omega\}$ be a sequence of countable covers of
   $\wb{Y}$. Write $Y$ as $\cup_{\alpha<\kappa}K_\alpha$ with $\kappa<\mathfrak{d}$ and $K_\alpha$ compact.
   For each $\alpha\in\kappa$,
   define $f_\alpha(n) = \min(\{ m\,:\, \cup_{\ell\le m}U_{k,\ell}\supset K_\alpha\,\forall k\le n\})$. 
   Since $\kappa<\mathfrak{d}$ there is some $f\in\,^\omega\omega$ such that for each $\alpha<\kappa$,
   $f(n)>f_\alpha(n)$ for infinitely many $n$.
   Hence $Y\subset \cup_n \cup_{k\le f(n)} U_{n,k}$, and $\wb{Y}$ is therefore weakly Menger.
   Since it is also od-Menger,
   Lemma \ref{lemma:DeltaOandOD} implies that $\wb{Y}$ is Menger.
\end{proof}

\begin{lemma}\label{lemma:unionb}
   Let $X$ be Lindel\"of and od-Menger, and $Y\subset X$ 
   be the union of $<\mathfrak{b}$ Menger subsets.
   Then $\wb{Y}$ is Menger. 
\end{lemma}
Classical argument again.
\begin{proof}
   Let $Y$ be $\cup_{\alpha <\kappa}Y$, where each $Y_\alpha$ is Menger and $\kappa<\mathfrak{b}$.
   Then $\wb{Y}$ is Lindel\"of.
   Let $\mathcal{U}_n = \{U_{n,m}\,:\,m\in\omega\}$ be a sequence of countable covers of it.
   For each $n\in\omega$ and each $\alpha <\kappa$,
   by Lemma \ref{lemma:classicalequiv} there is $f_\alpha(n)\in\omega$
   such that each $y\in Y_\alpha$ is in $\cup_{k\le f_\alpha(n)} U_{n,k}$ for infinitely many $n$.
   Since $\kappa<\mathfrak{b}$, there is some $f\in\,^\omega\omega$ such that,
   for each $\alpha<\kappa$, there is $m(\alpha)$ with $f(n)>f_\alpha(n)$ whenever $n\ge m(\alpha)$.
   We may assume that $f$ is strictly increasing.
   Set $\mathcal{F}_n = \{U_{n,m}\,:\,m\le f(n)\}$. Given $y\in Y_\alpha$,
   there is $n\ge m(\alpha)$ such that $y_\alpha\in \cup_{k\le f_\alpha(n)}U_{n,k}$.
   Hence, $y\in\cup \mathcal{F}_n$.
   This shows that $\cup_{n\in\omega}\cup\mathcal{F}_n$ is a cover of $Y$.
   It follows that $\wb{Y}$ satisfies $\ufin{\OO}{\DD}$ and is thus Menger by
   Lemma \ref{lemma:DeltaOandOD}.
\end{proof}

In a space $X$, let $\mathcal{M}_X$ be the ideal of subsets of $X$ whose closure is Menger
and $\mathcal{M}^w_X$ be the family of subsets of $X$ whose closure is weakly Menger.
$\mathcal{M}^w_X$ may fail to be an ideal, because closed subspaces of weakly Menger spaces
may fail to be weakly Menger or even weakly Lindel\"of, see e.g. the tangent disk space \cite[Ex. 82]{CEIT}. 
This cannot happen if $X$ is od-Menger.

\begin{cor} 
   \label{cor:ideals}
   If 
   $\ufin{\Delta_X}{\OO_X}$ holds, then  
   $\mathcal{M}^w_X=\mathcal{M}_X$ are $\sigma$-ideals, actually $\mathfrak{b}$-closed
   ideals.
\end{cor}
\proof
   If $E\subset X$ is such that $\wb{E}$ is weakly Menger, then by Lemma \ref{lemma:DeltaOandOD}
   $\wb{E}$ is Menger. Hence 
   $\mathcal{M}^w_X=\mathcal{M}_X$, and the rest follows by
   Lemma \ref{lemma:unionb}.
\endproof

\begin{q} Can $X$ be an increasing union of closed Menger subspaces if 
$\ufin{\Delta_X}{\OO_X}$ holds but $\ufin{\OO_X}{\OO_X}$ does not~?
Is it interesting to know the answer~?
\end{q}

Recall that given a cardinal $\kappa$, we say that
a space $X$ is {\em $\kappa$-concentrated on a subset $Y$} iff
$|X-U|<\kappa$ for each open $U\supset Y$. The following is classical
and, as far as we know, was first noted by Wingers \cite[Cor.3.9]{Wingers:1995}.
\begin{thm}\label{thm:Wingers}
   If $X$ is a Lindel\"of space $\mathfrak{d}$-concentrated on a Menger subspace $Y$,
   then $X$ is Menger as well.
\end{thm}

In order to have a concise statement of the next lemma, let us say until the end of this section
that a space is a {\em small union} iff it is either the union of $<\mathfrak{d}$ compact subsets,
or the union of $<\mathfrak{b}$ Menger subsets.
Say that a space $X$ is {\em smally-concentrated} 
on a subset $Y$ iff $X-U$ is a small union for each open $U\supset Y$.

\begin{lemma}
   Suppose that $X$ is smally concentrated on a subset $E$,
   and that $E\subset \wb{F}$, where $F$ is a small union.
   If $X$ is Lindel\"of and satisfies $\ufin{\Delta}{\OO}$ then it satisfies $\ufin{\OO}{\OO}$.
\end{lemma}
\proof
   $\wb{F}$ is Menger by Lemmas \ref{lemma:compact}--\ref{lemma:unionb}.
   Given a sequence of covers, we may thus cover $\wb{F}$ with finite families of the odd indices.
   The uncovered space is a closed small union and is thus 
   also Menger by Lemmas \ref{lemma:compact}--\ref{lemma:unionb}, 
   and we cover it with finite families of the even indices.
\endproof

Because singletons are undoubtedly compact, an immediate corollary is the following.
\begin{cor}
   If $X$ is a Lindel\"of space which satisfies $\ufin{\Delta}{\OO}$ but not $\ufin{\OO}{\OO}$
   and $X$ is $\mathfrak{d}$-concentrated on a subset $E$, then the density of $E$ is $\ge\mathfrak{d}$.
\end{cor}

A point in a space $X$ is a {\em $(\ge\kappa)$-accumulation point} of some $Y\subset X$
iff the intersection of each neighborhood of it with $Y$ has cardinality $\ge\kappa$.
If $X$ is a space, $Y\subset X$ and $\kappa$ is a cardinal, 
we denote by $\text{cl}_\kappa(Y,X)$ the subspace of $(\ge\kappa)$-accumulation points of $Y$.
One sees easily that $\text{cl}_\kappa(Y,X)$ is a closed subspace of $X$.
We denote $\text{cl}_\kappa(X,X)$ by $X_\kappa$.
A space is {\em linearly Lindel\"of} iff any chain cover has a countable subcover
or equivalently (as well known)
iff any subset of uncountable regular cardinality $\kappa$ has a 
$(\ge\kappa)$-accumulation point.

\begin{lemma}\label{lemma:kappa-concentrated}
   If $X$ is linearly Lindel\"of and $\kappa$ is regular and uncountable, then $X$ is $\kappa$-concentrated on $X_\kappa$.
\end{lemma}
\proof
   Let $U\supset X_\kappa$ be open. Then $X-U$ is closed and linearly Lindel\"of, hence, if its cardinality is 
   $\ge\kappa$ it contains 
   a point of $X_\kappa$, a contradiction.
\endproof

\begin{lemma}
   \label{lemma:X_bMenger}
   Let $X$ be Lindel\"of. Then the following hold.\\
   (1) 
   $ \ufin{\OO_{X_{\mathfrak{b}}}}{\OO_{X_{\mathfrak{b}}}}\,\Longrightarrow\,
      \ufin{\OO_X}{\OO_X}$.\\
   (2)  $ \ufin{\Delta_X}{\OO_X}\text{ and } 
      \ufin{\OO_{X_{\mathfrak{b}}}}{\DD_{X_{\mathfrak{b}}}} \,\Longrightarrow\,
      \ufin{\OO_X}{\OO_X}.
      $
\end{lemma}
\proof
   \ \\
   (1) 
   By Theorem \ref{thm:Wingers}, Lemma \ref{lemma:kappa-concentrated} 
   and the fact that $\mathfrak{b}$ is regular and $\le\mathfrak{d}$.\\
   (2) Since $X_{\mathfrak{b}}$ is closed, $\ufin{\Delta_{X_{\mathfrak{b}}}}{\OO_{X_{\mathfrak{b}}}}$ holds.
   Hence,
   so does $\ufin{\OO_{X_{\mathfrak{b}}}}{\OO_{X_{\mathfrak{b}}}}$ by Lemma \ref{lemma:DeltaOandOD},
   and we finish with (1).
\endproof

\begin{q} Is it useful to iterate (\`a la Cantor) the procedure
          $X_\mathfrak{b}$, $\left(X_\mathfrak{b}\right)_\mathfrak{b}$, 
          $\left(\left(X_\mathfrak{b}\right)_\mathfrak{b}\right)_\mathfrak{b}$, \dots~? 
          What happens (to the Menger-like properties) after a countable number of steps~?
\end{q}

We now summarize the results of this section in one Theorem.
Point (1) was not discussed before but is easy, so we added it as well.
\begin{thm} 
   \label{thm:summary}
   Let $X$ be Lindel\"of, od-Menger and non-Menger,
   and $\omega+1$ be the convergent sequence space. Then the following hold.\\
   \indent (1) $X\times (\omega+1)$ is non-od-Menger and the union of two od-Menger subspaces.\\
   \indent (2) $X_{\mathfrak{b}}$ is not weakly Menger.\\
   If $X$ is moreover everywhere non-Menger, then the following hold as well.\\
   \indent (3) A Menger subspace of $X$ is nowhere dense and has Menger closure.\\
   \indent (4) A closed subspace of $X$ is Menger iff it is nowhere dense.\\
   \indent (5) The union of $<\mathfrak{b}$ Menger subspaces of $X$ is nowhere dense and has Menger closure.\\
   \indent (6) The union of $<\mathfrak{d}$ compact subspaces of $X$ is nowhere dense and has Menger closure.\\
   \indent (7) $X$ is Baire, actually, the intersection of $<\mathfrak{b}$ open dense sets is dense.
\end{thm}
(We are aware that (5)$\Rightarrow$(3) trivially, but we prefer to separate the statements.)
\proof
  (1) $X\times\{\omega\}$ is non-Menger, closed and nowhere dense in $X\times (\omega+1)$, which is therefore not od-Menger. 
      It is easy to see that
      $X\times \omega$ is od-Menger since $\omega$
      is discrete, and $X\times (\omega+1) = X\times \omega\sqcup X\times\{\omega\}$. \\
  (2) Immediate by Lemma \ref{lemma:X_bMenger} (2).\\
  (3) Follows by (5), but a direct proof is almost immediate. Indeed, let $E\subset X$ be Menger.
      By Lemma \ref{lemma:DeltaOandOD}, $\wb{E}$ is Menger.
      Since $M(X)=\varnothing$, $\wb{\text{int}(\wb{E})} = \varnothing$, hence
      $E$ has empty interior.\\
  (4) Immediate by definition and the fact that $X$ is od-Menger.\\
  (5) By Corollary \ref{cor:ideals} and (3).\\
  (6) Immediate by Lemma \ref{lemma:compact} and (3).\\
  (7) Follows directly from (4) and (5).
\endproof


\section{Examples of $\ufin{\Delta}{\OO} \not\Longrightarrow \ufin{\OO}{\OO}$}

While answering a question on mathoverflow, the second author showed that 
the space $\mathcal{K}(\PP)$ 
($\PP$ the irrationals), which is a small variation of a space defined by 
E. van Douwen, F. Tall and W. Weiss in the 70s in \cite[Theorem 3]{vanDouwenTallWeiss:1977}, 
does not satisfy $\ufin{\OO}{\DD}$. 
He also observed that this space has a dense Luzin subspace $\mathcal{L}$ under {\bf CH},
which yields an example of an od-Menger non-Menger space.
As said in the introduction,
another example is implicit in a paper of M. Sakai in \cite[Thm 2.4 (b)]{Sakai:2014},
but does not have a $G_\delta$ diagonal, while $\mathcal{L}$ has
(enabling us to prove Corollary \ref{ex:main} below).
We now give the details of the definitions and proofs of this latter example.\\
If $X$ is a space,
the space $\mathcal{K}(X)$ is defined as the set of compact nowhere dense sets of $X$,
with the Pixley-Roy topology. That is,
for any $F,G\subset X$, set
$$ [F,G] = \{ H\in \mathcal{K}(X)\,:\, F\subset H\subset G\};$$
a base of neighborhoods of $K\in\mathcal{K}(X)$
is given by
$$ \{ [K,U]\,:\,K\subset U,\,U\subset X\text{ is open}\}. $$ 
The space defined in \cite[Theorem 3]{vanDouwenTallWeiss:1977} is
$\mathcal{K}(\mathbb{R})$, where it is shown that it
is first countable, ccc, $0$-dimensional, has a $G_\delta$-diagonal
and a $\sigma$-centered base (see the reference for definitions).
We now show that the same holds for its subspace $\mathcal{K}(\PP)$.
First, a definition:
we say that a space $X$ is {\em DTW-complete} iff there is a 
sequence $\{\mathfrak{B}_n\,:\,n\in\omega\}$ of families of closed subsets of $X$
   such that
   \begin{itemize}
     \item[(DTW 1)] 
     for each $n\in\omega$, for each $K\in X$ and for each neighbourhood $\mathcal{U}$
     of $K$ in $X$, there is a $\mathcal{B}\in\mathfrak{B}_n$ such that
     $\displaystyle K\in \text{int}_{X}(\mathcal{B})\subset \mathcal{B}\subset\mathcal{U}$; 

     \item[(DTW 2)] if $\mathfrak{U}\subset\cup\{\mathfrak{B}_n\,:\,n\in\omega\}$ 
     is any centered family 
     such that $\mathfrak{U}\cap\mathfrak{B}_n\not=\varnothing$ for each $n\in\omega$, then
     $\cap \mathfrak{U}\not=\varnothing$.
   \end{itemize}
(Recall that a family is centered iff any finite subfamily has a non-empty intersection.)
The proof that a DTW-complete space is Baire is so straighforward that
it was not given in \cite{vanDouwenTallWeiss:1977}.
Let us give a version of this proof in a language appropriate for 
generalizations under Martin's axiom (we note 
that this generalization is also well known, although we did not find a proof in print,
and we use it in the last remark of this paper).
\begin{lemma}
   \label{lemma:KPBaire}
   A DTW-complete space is Baire.
\end{lemma}
\proof
   Let $\mathsf{P}$ be the poset $\cup_n\mathfrak{B}_n$ (where $\{\mathfrak{B}_n\,:\,n\in\omega\}$ is 
   as in the definition of DTW-completeness), 
   ordered as $\mathcal{B}_0\ge\mathcal{B}_1$ iff $\text{int}(\mathcal{B}_0)\subset\text{int}(\mathcal{B}_1)$.
   Given a countable collection of open dense sets $\mathcal{D}_n$,
   set $\mathsf{D}_n = \{ \mathcal{B}\in\mathsf{P}\,:\,\mathcal{B}\subset \mathcal{D}_n\}$,
   then (DTW 1) implies that each $\mathsf{D}_n$ is dense in $\mathsf{P}$.
   By the Rasiowa-Sikorski Lemma, there is a filter $\mathsf{G}$ intersecting each $\mathcal{D}_n$.
   Then (DTW 2) ensures that
   $\cap \mathsf{G}\not=\varnothing$. 
   Any point in this intersection is a member of $\cap_{n\in\omega}\mathcal{D}_n$.
\endproof

\begin{lemma}
  \label{lemma:KPproperties}
  $\mathcal{K}(\PP)$ is first countable, ccc, Baire, $0$-dimensional, of cardinality $\mathfrak{c}$, 
  has a $G_\delta$-diagonal
  and a $\sigma$-centered base, and has no isolated points.
\end{lemma}
\begin{proof}
The fact that no point of $\mathcal{K}(\PP)$ is isolated is clear.
Notice that $\mathcal{K}(\PP)\subset \mathcal{K}(\mathbb{R})$: a compact subspace of $\PP$ is also compact in $\R$.
Moreover, by definition the topology of $\mathcal{K}(\PP)$ is equal to the subspace topology.
First countability, $0$-dimensionality and having a $G_\delta$ diagonal are hereditary,
and the proof in \cite[Theorem 3]{vanDouwenTallWeiss:1977} 
of the fact that $\mathcal{K}(\mathbb{R})$ has a $\sigma$-centered base (which implies ccc)
only uses the fact that $\R$ is second countable and thus also works for $\mathcal{K}(\PP)$.
Since a second countable space has at most $\mathfrak{c}$ closed sets, $\mathcal{K}(\PP)$ has cardinality $\mathfrak{c}$.
Van Douwen, Tall and Weiss also show in \cite[Theorem 3]{vanDouwenTallWeiss:1977} that
$\mathcal{K}(\R)$ is DTW-complete. We repeat their proof (with a very small adaptation)
and show that $\mathcal{K}(\PP)$ is also DTW-complete, and hence Baire.
For this, let $\{B_n\,:\,n\in\omega\}$ be a base for $\mathbb{R}$ and $\{q_n\,:\,n\in\omega\}$ be an enumeration
of $\mathbb{Q}$.
We set 
\begin{align*}
   \mathfrak{B}_n = \{ [K,F\cap\PP]\,:\, 
   & K\in\mathcal{K}(\PP),\,F\subset\mathbb{R}\text{ is compact in }\mathbb{R},\\
   & K\subset\text{int}_{\mathbb{R}}(F),\, B_n-F\not=\varnothing,\,q_n\not\in F\}.
\end{align*}
Then (DTW 1) is immediate, since any $K\in\mathcal{K}(\PP)$ is compact and nowhere dense and
does obviously not contain $q_n$.
Let now $\mathfrak{U}\subset\cup\{\mathfrak{B}_n\,:\,n\in\omega\}$ be a centered family 
such that $\mathfrak{U}\cap\mathfrak{B}_n\not=\varnothing$ or each $n\in\omega$.
Thus $\mathfrak{U} = \{[K_i,F_i\cap\PP]\,:\,i\in I\}$. Set $F = \cap_{i\in I} F_i$ (seen 
as a subset of $\mathbb{R}$). Then for each $n\in\omega$, $B_n-F\not=\varnothing$
and $q_n\not\in F$, since
there is some $i\in I$ such that $B_n-F_i\not=\varnothing$ and $q_n\not\in F_i$. 
Thus, $F$ is a nowhere dense compact subset of $\PP$.
Now, since $\mathfrak{U}$ is centered, for any $i,j\in I$ we have
$$ 
  \varnothing \not= [K_i,F_i\cap\PP] \cap [K_j,F_j\cap\PP] = [K_i\cup K_j,F_i\cap F_j\cap\PP],
$$
hence $K_i\subset F_j\cap\PP$ for any $i,j\in I$. It follows that  
$K_i\subset F\subset F_i\cap\PP$ for each $i$,
which implies that $\varnothing\not= F\in \cap\mathfrak{U}$ and (DTW 2) holds.
This finishes the proof of Theorem \ref{lemma:KPproperties}.
\end{proof}

If we could add the word ``dense'' in the lemma below,
we would have a shortcut for the proof that $\mathcal{K}(\PP)$ is Baire, but unfortunately it is not the case.
\begin{lemma}\label{lemma:KGdelta}
   $\mathcal{K}(\PP)$ is a $G_\delta$ subset of $\mathcal{K}(\mathbb{R})$.
\end{lemma}
\proof
   $\mathcal{E}_p = \{K\in\mathcal{K}(\R)\,:\,K\not\ni p\}$ is open in $\mathcal{K}(\R)$ for each $p\in\R$,
   since for each $K\not\ni p$,  $[K,\R-\{p\}]\subset \mathcal{E}_p$.
   Since each $K\in \mathcal{K}(\PP)$ is compact in $\PP$,
   it is closed in $\R$ as well. It follows that
   $\mathcal{K}(\PP)=\cap_{p\in\Q} \mathcal{E}_p$ is a $G_\delta$.  
\endproof 

Recall that a {\em Luzin space} is an uncountable space without isolated points such that 
any closed nowhere dense set is countable; $\sone{\Delta}{\OO}$ holds immediately
for Luzin spaces. Luzin spaces are hereditarily Lindel\"of (see e.g. \cite[p. 308]{Roitman:1984}).

\begin{thm}
  Assume {\bf CH}.
  Then there is a hereditarily Lindel\"of first countable $0$-dimensional space without isolated points and
  with a $G_\delta$ diagonal
  satisfying $\sone{\Delta}{\OO}$ but not $\ufin{\OO}{\DD}$.\\  
  \label{thm:spadaro}
\end{thm}
\proof
By \cite[Corollary to Theorem 1]{vanDouwenTallWeiss:1977}, under {\bf CH}
every uncountable ccc first countable Baire 
space without isolated points 
possesses a dense Luzin subspace.
The proof is by induction of length $\omega_1$, using Baire's property for the intermediate stages.
Let thus $\mathcal{L}$ be a dense Luzin subspace of $\mathcal{K}(\PP)$.
The fact that $\sone{\Delta}{\OO}$ holds is immediate, $0$-dimensionality, first countability 
and having a $G_\delta$
diagonal are hereditary, and recall that 
Luzin spaces are hereditarily Lindel\"of.
Since $\mathcal{K}(\PP)$ has no isolated points and $\mathcal{L}$ is dense, $\mathcal{L}$ has no isolated points either.
The only remaining thing to show is that $\ufin{\OO_{\mathcal{L}}}{\DD_{\mathcal{L}}}$
does not hold. 
\\
Since $\mathcal{L}$ is dense in $\mathcal{K}(\PP)$, it is enough to show that 
 the later is not weakly Menger.
 Since $\PP$ is not Menger, 
 there is a countable sequence $\{\mathcal{U}_n\,:\,n\in\omega\}$
 of open covers of $\PP$ which witnesses that. 
 For every $K\in\mathcal{K}(\PP)$
 let $\mathcal{U}_n^K$ 
 be a finite subcollection of $\mathcal{U}_n$ 
 which covers $K$. 
 Then $\mathfrak{O}_n= \{ [K , \cup \mathcal{U}_n^K]\,:\, K\in \mathcal{K}(\PP)\}$
 is an open cover of $\mathcal{K}(\PP)$ for every $n\in\omega$.
 Let $\mathfrak{G}_n$ be a finite subcollection of $\mathfrak{O}_n$ for each $n$,
 and let $\mathfrak{F}_n$ be the finite subset of $\mathcal{K}(\PP)$ such that 
 $K\in\mathfrak{F}_n$ iff $[K , \cup \mathcal{U}_n^K]\in\mathfrak{G}_n$.
 Then $\cup \{\mathcal{U}_n^K\,:\,K\in \mathfrak{F}_n\}$ is a finite subcollection of 
 $\mathcal{U}_n$ for every $n\in\omega$, and therefore there is a point
 $y\in\PP - \left(\cup\left\{\cup \mathcal{U}_n^K\,:\,K\in \mathfrak{F}_n\right\}\right)$.
 It follows that $[\{y\},\PP]$ is a non-empty open subset of $\mathcal{K}(\PP)$
 which is disjoint from $\cup\{\cup\mathfrak{G}_n\,:\,n\in\omega\}$, 
 thus showing that $\mathcal{K}(\PP)$ is not weakly Menger. 
\endproof

As above, we let $\omega+1$ be the convergent sequence space.
\begin{cor}[{\bf CH}]
   \label{ex:main}
   There is a refinement of the Euclidean topology on a subset $E$ of the reals such that 
   $E$ is hereditarily Lindel\"of, $0$-dimensional, first countable and:\\
   (1) $E$ satisfies $\sone{\Delta}{\OO}$ but not $\ufin{\OO}{\DD}$;\\
   (2) $E\times (\omega+1)$ is Lindel\"of and the union of
       $2$ od-Menger spaces, but is not od-Menger.
\end{cor}
\begin{proof}
   Theorem \ref{thm:spadaro} and
   Proposition \ref{prop:real_example} provide $E$, the rest is shown in Theorem \ref{thm:summary}.
\end{proof}

Recall that the product of a Menger space with a compact space is Menger
(to our knowledge, this was noticed first by Telg\'arski in 1972 \cite[Lemma 3.6]{Telgarsky:1972}). 
The previous corollary shows that this does not hold
(consistently) for od-Menger spaces. The fact that a Lindel\"of space which is the union of two od-Menger
spaces may not be od-Menger shows a stark difference with Menger spaces.

\begin{q} 
   Is there a Lindel\"of od-Menger non-Menger space in {\bf ZFC}~?
\end{q}

We now briefly discuss the example of Lindel\"of od-Menger non-Menger space
due to 
M. Sakai in \cite[Thm 2.4 (b)]{Sakai:2014} alluded to in the introduction. 
Sakai's space, presented as an example of a Lindel\"of non-weakly Menger space, 
also needs {\bf CH} and
is similar to that of Theorem \ref{thm:spadaro} (that is: a Luzin non-weakly Menger space,
which therefore satisfies $\sone{\Delta}{\OO}$).
Sakai actually does not state that his space is Luzin, because he is only interested in its
hereditary Lindel\"ofness, but its existence is a consequence of his Lemma 2.3,
which, although quoted from another source, is exactly
\cite[Theorem 1]{vanDouwenTallWeiss:1977}, and the latter yields Luzin spaces.
We are sure that if \cite{meszigues-od-sel} had been available to 
him when he wrote his paper, Sakai 
would have noticed that his space gives an example of Lindel\"of od-Menger non-Menger space.

Sakai actually presents another Lindel\"of non-weakly Menger space in his paper, which is much simpler
\cite[Thm 2.5]{Sakai:2014} and, unlike the first one, does have a $G_\delta$ diagonal.
However, this space does not satisfy $\ufin{\Delta}{\OO}$,
let us show why. 
\\
Let $X$ be a space, $\tau$ its topology, and let $Y\subset X$. 
Set $X_Y$ to be the topological space $X$ with the topology generated by 
$(X-Y)\cup\tau$. Thus an open set of $X_Y$ is the union of a subset of $X-Y$ and $U\cap Y$ for some $U\in\tau$
and the induced topology on $Y\subset X_Y$ is the same as that of $X$.
Sakai's second example 
is $K_B$ where $K$ is the Cantor space and 
$B$ is a Bernstein set. The following simple lemma shows that no Lindel\"of 
such space 
satisfies $\ufin{\Delta}{\OO}$ but not
$\ufin{\OO}{\OO}$.

\begin{lemma}\label{lemma:mitchell_style}
   Let $Y\subset X$ be spaces such $\ufin{\OO_Y}{\DD_Y}$ holds.
   If $X_Y$ is Lindel\"of and satisfies $\ufin{\Delta}{\OO}$, then it also satisfies $\ufin{\OO}{\OO}$.
\end{lemma}
\proof
   Since $X_Y$ is Lindel\"of, $X$ is $\aleph_1$-concentrated on $Y$: if $U$ is
   any open set containing $Y$, then $|X-U|<\aleph_1$ (because it is closed discrete). 
   Since $X_Y$ satisfies $\ufin{\Delta}{\OO}$, so does its closed subset $Y$, and by Lemma \ref{lemma:DeltaOandOD}
   $\ufin{\OO_Y}{\OO_Y}$ holds. We conclude with Theorem \ref{thm:Wingers}.
\endproof

\begin{thm}[{M. Sakai \cite[Thm 2.5]{Sakai:2014}}] 
\label{thm:Sakai2}
There is a regular Lindel\"of space with a $G_\delta$-diagonal 
which satisfies neither $\ufin{\OO}{\DD}$ nor $\ufin{\Delta}{\OO}$.
\end{thm}
\proof
Sakai's example is $K_B$ where $K$ is the Cantor space and $B$ is a Bernstein set in $K$. 
It is shown in \cite[Thm 2.5]{Sakai:2014} that $K_B$ is not weakly Menger.
Since $B$ is separable, $\ufin{\OO_B}{\DD_B}$ holds.
Lemma \ref{lemma:mitchell_style}
implies that $K_B$ cannot
satisfy $\ufin{\Delta}{\OO}$.
\endproof

\vskip .3cm
To finish, we notice that under Martin's axiom {\bf MA},
the space $\mathcal{K}(\PP)$ has the property of being {\em strongly Baire}, that is:
the intersection of $<\mathfrak{c}$ open dense subsets is dense (this is essentially Lemma \ref{lemma:KPBaire},
together with the fact that $\mathcal{K}(\PP)$ is ccc).
We say that a space is {\em $\mathfrak{c}$-Luzin} iff it has cardinality $\mathfrak{c}$ while
its closed nowhere dense subsets are of cardinality $<\mathfrak{c}$.
It was noted a long time ago (see e.g. the remarks between pages 206 and 207 in \cite{Miller:1984})
that 
the proof of \cite[Theorem 1]{vanDouwenTallWeiss:1977} under {\bf CH} generalizes 
almost automatically
to yield the following result. 

\begin{thm}
   Let $X$ be a ccc strongly Baire space of cardinality $\mathfrak{c}$ with a $\pi$-base of cardinality $\le\mathfrak{c}$.
   Then $X$ contains a dense $\mathfrak{c}$-Luzin subset.
   \label{thm:vDTWMA}
\end{thm}
Hence, under {\bf MA}, we obtain a dense $\mathfrak{c}$-Lusin subspace $\mathcal{L}$
of $\mathcal{K}(\PP)$.
Since any nowhere dense closed subset of $\mathcal{L}$ has cardinality $<\mathfrak{c}$
and $\mathfrak{d}=\mathfrak{c}$ under {\bf MA}, $\mathcal{L}$ would be od-Menger by 
Theorem \ref{thm:led}
{\em if} it were Lindel\"of.
Unfortunately, we do not have any indication that it is the case.

{\footnotesize
\vskip .4cm
\noindent
Mathieu Baillif \\
Haute \'ecole du paysage, d'ing\'enierie et d'architecture (HEPIA) \\
Ing\'enierie des technologies de l'information TIC\\
Gen\`eve -- Suisse\\
Email: {\tt mathieu.baillif@hesge.ch}
\
\vskip .3cm \noindent
Santi Spadaro\\
Department of Mathematics and Computer Science\\
University of Catania\\
viale Andrea Doria n. 6\\
95125 Catania, Italy
\\
Email: {\tt santi.spadaro@unict.it}
}

\end{document}